\documentclass[conference]{ieeeconf}

\renewcommand{\textfloatsep}{5mm} 

\IEEEoverridecommandlockouts

\usepackage[acronym]{glossaries}
\usepackage{amsmath,amssymb,amsfonts}
\usepackage{algorithmic}
\usepackage[pdftex]{graphicx}
\usepackage{xcolor}
\def\BibTeX{{\rm B\kern-.05em{\sc i\kern-.025em b}\kern-.08em
    T\kern-.1667em\lower.7ex\hbox{E}\kern-.125emX}}

\usepackage{units}
\usepackage{mathtools} 
\usepackage{amsthm}
\usepackage{lettrine}

\usepackage{scrextend} 

\usepackage{tabularx}
\usepackage[english]{babel}
\usepackage{cite}

\usepackage{array}
\usepackage{color}
\usepackage{amssymb}
\usepackage{multicol}

\DeclareUnicodeCharacter{2212}{-}

\usepackage{todonotes}

\usepackage{tikz}
\usetikzlibrary{arrows,calc}
\usepackage{pgfplots}
\usepgfplotslibrary{groupplots}

\usepackage{amsthm}
\usepackage{algorithm}
\usepackage{hyperref}
\usepackage{epstopdf}
\usepackage{mathrsfs}  

\newcounter{thm}
\newtheoremstyle{mystyle}
  {}
  {}
  {}
  {}
  {\bfseries\color{black}}
  {.}
  {\newline}
  {\thmname{#1}\thmnumber{ #2: }\normalfont\color{black}\thmnote{ \textit{(#3)}}}
\theoremstyle{mystyle}
\newtheorem{prob}[thm]{Problem}

\newtheorem{theorem}{Theorem}[section]

\newacronym{acr:cvt}{CVT}{continuously variable transmission}
\newacronym{acr:dp}{DP}{dynamic programming}
\newacronym{acr:ecms}{ECMS}{equivalent consumption minimization strategies}
\newacronym{acr:eltms}{ELTMS}{equivalent lap time minimization strategies}
\newacronym{acr:em}{EM}{electric motor}
\newacronym{acr:es2k}{ES2K}{Energy Storage to Kinetic}
\newacronym{acr:F1}{F1}{Formula 1}
\newacronym{acr:FIA}{FIA}{F\'{e}d\'{e}ration Internationale de l'Automobile}
\newacronym{acr:fgt}{FGT}{fixed-gear transmission}
\newacronym{acr:FD}{FD}{final drive}
\newacronym{acr:ice}{ICE}{internal combustion engine}
\newacronym{acr:k2es}{K2ES}{Kinetic to Energy Storage}
\newacronym{acr:mgu}{MGU}{motor generator unit}
\newacronym{acr:mguh}{MGU-H}{motor generator unit heat}
\newacronym{acr:mguk}{MGU-K}{motor generator unit kinetic}
\newacronym{acr:mpc}{MPC}{model predictive control}
\newacronym[description={Energy Management Strategy}, \glslongpluralkey={Energy Management Strategies},\glsshortpluralkey={EMSs}]{EMS}{EMS}{Energy Management Strategy}%
\newacronym{acr:pmp}{PMP}{Pontryagin's Minimum Principle}
\newacronym{acr:pu}{PU}{power unit}
\newacronym[description={Powertrain Operation}, \glslongpluralkey={Powertrain Operations},\glsshortpluralkey={POs}]{acr:PO}{PO}{Powertrain Operation}%
\newacronym{acr:rmse}{RMSE}{root-mean-square error}
\newacronym{acr:socp}{SOCP}{second-order cone program}
\newacronym{acr:soe}{SoE}{State of Energy}

\newcommand{\pushright}[1]{\ifmeasuring@#1\else\omit\hfill$\displaystyle#1$\fi\ignorespaces}
\newcommand{\pushleft}[1]{\ifmeasuring@#1\else\omit$\displaystyle#1$\hfill\fi\ignorespaces}
\makeatother
\newif\ifmargincomments 
\margincommentstrue

\newif\ifextendedversion 
\extendedversionfalse

\ifmargincomments

\else

\fi
\begin{document}

\title{\bf A Computationally Efficient and Human Implementable\\ Minimum-lap-time Control Policy for Energy-limited Race Cars}
			

\author{Erik van den Eshof, Wytze de Vries, Mauro Salazar%
\thanks{
	The authors are with the Control Systems Technology section, Department of Mechanical Engineering, Eindhoven University of Technology (TU/e), Eindhoven, 5600 MB, The Netherlands.
E-mails: {\tt\footnotesize \{r.c.p.v.d.eshof,w.a.b.d.vries,m.r.u.salazar\}@tue.nl}
}}

\maketitle
\thispagestyle{plain}
\pagestyle{plain}
\begin{abstract}
	This paper presents a provably optimal, real‑time capable energy management policy for race cars that provides simple human-driver-implementable control cues. Specifically, we first formulate the energy-constrained minimum-lap-time control problem via Pontryagin's Minimum Principle (PMP) and derive the optimal policy and costate dynamics using Karush–Kuhn–Tucker (KKT) optimality conditions.
	We show that the optimal control policy follows a bang-bang structure that is easily implementable by a human driver, eliminating the need for potentially dangerous active throttle pedal overwrites or distracting signals.
	Moreover, the analytical formulation of the optimal system dynamics allows us to recast the problem as a sequence of boundary‑value problems, which can be efficiently solved using root‑finding methods.
	Our results show that our proposed approach can compute the same globally optimal control strategies of existing numerical methods based on direct optimal control, whilst drastically reducing computation time from the order of seconds to milliseconds.
\end{abstract}

\section{Introduction}
The importance of energy management in motorsports is more significant than ever with the introduction of in-race charging in Formula E in 2025 and new 2026 Formula~1 regulations featuring a heavy reliance on battery energy deployment~\cite{F12026TechRegs,FE20252026TechRegs}. The high-speed nature of racing calls for computationally efficient energy management algorithms that operate safely with human drivers, minimizing control overrides or distractions. Simultaneously, given the characteristically narrow performance margins in motorsports, the use of provably optimal methods is of paramount importance.\\
In situations where energy limitations are minor, optimal energy management strategies involve only small adjustments in deployed power, leading to small and predictable variations in vehicle speed. For safety reasons, however, active overwrites of the driver's power demand are often subject to strict rate limits or are prohibited entirely, ensuring that the driver remains in control of the vehicle. This constraint becomes limiting when energy availability is tight, such as in fully electric powertrains, where a race‑time‑optimal strategy may require substantial power reductions, and thus noticeable speed changes, to conserve energy~\cite{KampenHerrmannEtAl2023}. To maintain safe operation when driving in close proximity to competitors, the driver must stay in control with minimal cognitive load. Unlike tire grip limits, which drivers can physically perceive, energy limits cannot be sensed. Consequently, in races with pronounced battery or fuel limitations, drivers rely on externally provided cues or strategies to manage energy usage: the energy management system instructs the driver when to lift off the throttle or enable regenerative braking, allowing the driver to execute these actions only when it is safe to do so.

\setlength{\fboxsep}{0pt}
\setlength{\fboxrule}{1pt}
\begin{figure}[t]
	\centering
	\vspace{2mm} 
	\framebox{\includegraphics[width=\linewidth]{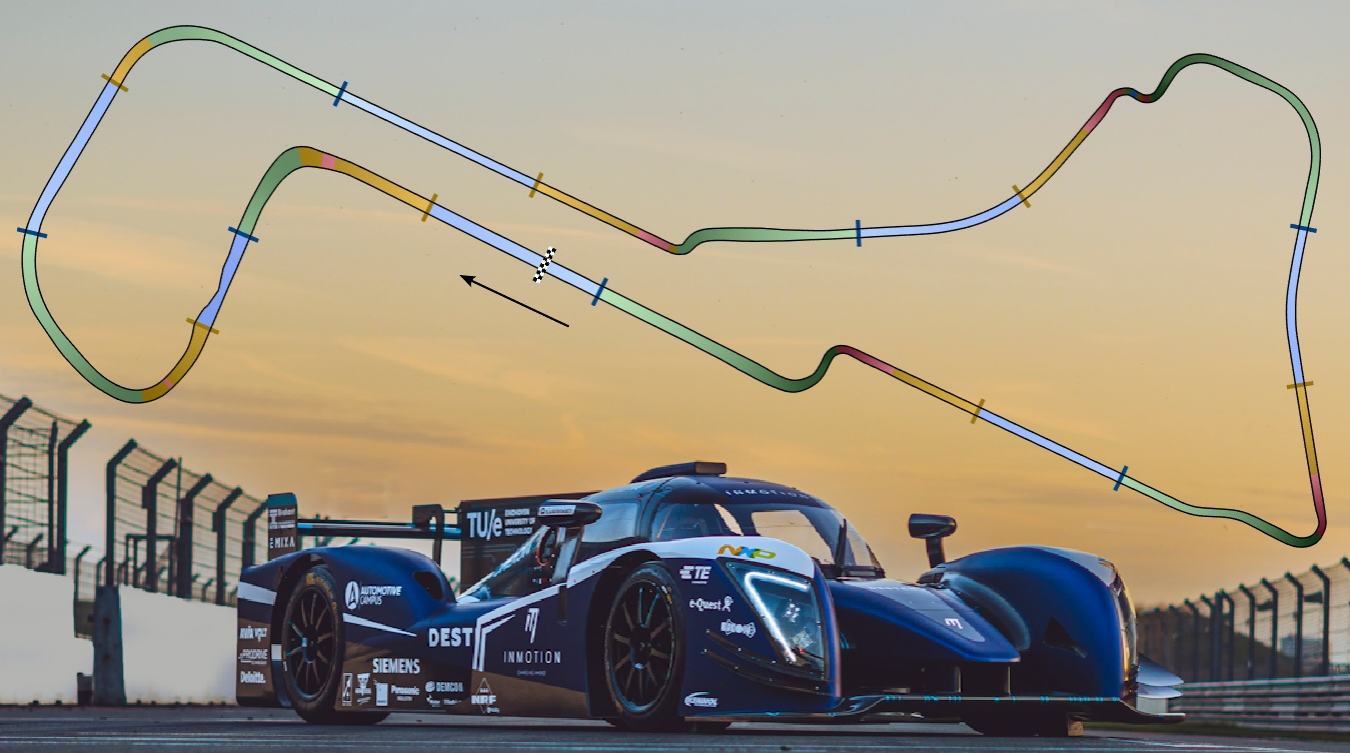}}
	\caption{Electric endurance racing car (InMotion)~\cite{InMotion}, typically severely energy-limited. Optimal coasting and regenerative braking points along the circuit are highlighted in blue and gold respectively. The dark and light regions respectively indicate slow and high speeds.}
	\label{fig:track}
	\vspace{-3mm}
\end{figure}
\paragraph*{Related Literature}
Energy management has become an increasingly important research area in motorsport as reliance on battery‑electric energy grows. Unlike in road vehicles, where the primary goal is to minimize energy consumption, motorsport focuses on minimizing lap-time under strict energy constraints~\cite{BorsboomSalazarEtAl2024}. Researchers have studied energy-limited racing applications of techniques such as convex optimization~\cite{EbbesenSalazarEtAl2018,KampenHerrmannEtAl2023,BoydVandenberghe2004,BorsboomFahdzyanaEtAl2021}, using simplified convex models, and nonlinear-programming (NLP)~\cite{LiuFotouhiEtAl2020,LimebeerPerantoniEtAl2014,DuhrBuccheriEtAl2023}, capturing more advanced dynamics, but lacking computational efficiency for real-time applications. For real-time optimal control, researchers applied sequential convex programming methods (SCP)~\cite{HerrmannSauerbeckEtAl2021,EshofKampenEtAl2025b}, and convex models to derive optimal policies using Pontryagin's minimum principle (PMP)~\cite{SalazarElbertEtAl2017}, with disturbance rejection through model-predictive control (MPC)~\cite{SalazarBalernaEtAl2017} and equivalent lap-time minimization strategies (ELTMS)~\cite{SalazarBalernaEtAl2018}, or a combination of both~\cite{NeumannFieniEtAl2023}. These methods are mostly limited to hybrid powertrains, where the power split between electric and combustion engine is actively controlled and the driver remains in control of most of the power delivered to the wheels. For more significant energy limitations, as is the case for energy management in races with refueling or (particularly) recharging stops, optimal control requires an active overwrite of the power delivered to the wheels~\cite{HerrmannSauerbeckEtAl2021,KampenHerrmannEtAl2023}. While suitable for autonomous racing applications, this is typically restricted by racing regulations and difficult to implement in a safe manner when there are human drivers involved~\cite{F12026TechRegs,FE20252026TechRegs}. In our previous work~\cite{EshofKampenEtAl2025}, we proposed to solve this problem by adapting a smooth, convex optimal control solution to generate close-to-optimal bang-bang control inputs to be executed by the driver, whereby the driver remains in control of the power delivery. However, this method lacks global optimality and only considers a coasting phase for energy saving, meaning regenerative braking is only used in combination with friction brakes, limiting energy recovery potential.\\
To the best of the author's knowledge, there are no real-time capable methods for provably optimal energy management of racing cars, that are executable by a human driver.
\paragraph*{Statement of Contributions}
This paper presents a provably optimal, real-time capable method for energy management of racing cars. Through an analytical derivation of the Karush–Kuhn–Tucker (KKT) conditions~\cite{PapalambrosWilde2017}, which are sufficient for global optimality of the convex energy management problem, we gain deeper insight into the behavior of the optimal solution. Hereby we find that by selecting linearized powertrain loss models, the optimal solution exhibits a desirable bang-bang control policy as a function of the kinetic energy to battery energy costate ratio, which can be followed by human drivers with little distraction and avoids the need for a potentially unsafe active throttle request overwrite. Subsequently, we apply constrained PMP to derive the optimal costate dynamics~\cite{Bertsekas1995,KeulenGillotEtAl2014} and frame the problem as a set of sequential boundary value problems, which we can efficiently solve using root-finding methods.

\section{Methodology}
In this section, we derive the minimum-lap-time control policy and optimal dynamics for the energy limited racing problem using the KKT conditions and PMP. First, we define our optimization problem.\\
In previous work we showed that energy management strategies have little influence on the driven spatial trajectory~\cite{EshofKampenEtAl2025b}. Therefore, we assume a fixed, pre-computed curvature-based trajectory in the formulation of our energy management problem, significantly simplifying the control problem by exclusively focusing on longitudinal dynamics, without the consideration of spatial trajectory optimization. We operate in the distance domain of a single racing lap $s\in[s_\mathrm{0},s_\mathrm{f}]$, and define our objective lap-time $t_\mathrm{lap}$ as the inverse of velocity $v$ integrated over this distance:
\par\nobreak\vspace{-5pt}
\begingroup
\allowdisplaybreaks
\begin{small}
	\label{eqn1}
	\begin{equation}
		t_\mathrm{lap} = \int_{s_\mathrm{0}}^{s_\mathrm{f}}\frac{\mathrm{d}t}{\mathrm{d}s}(s)\, \mathrm{d}s = \int_{s_\mathrm{0}}^{s_\mathrm{f}}\frac{1}{v(s)}\, \mathrm{d}s = \int_{s_\mathrm{0}}^{s_\mathrm{f}}\frac{1}{\sqrt{2\,E_\mathrm{kin}(s)/m}}\, \mathrm{d}s,\\
	\end{equation}
	
\end{small}
\endgroup
\noindent where $E_\mathrm{kin}$ is the kinetic energy and $m$ is the vehicle mass. We select kinetic energy as a state, as this results in favorable linear force-dynamics in the distance domain. Accordingly, our state and dynamics are defined as follows:
\par\nobreak\vspace{-5pt}
\begingroup
\allowdisplaybreaks
\begin{small}
	
	\begin{equation}
			\label{eqn2}
		\textstyle\frac{\mathrm{d}E_\mathrm{kin}}{\mathrm{d}s}(s) = F_\mathrm{m}(s) - F_\mathrm{brk}(s) - F_\mathrm{d}(s,\kappa,E_\mathrm{kin}),
	\end{equation}
	\begin{equation}
			\label{eqn3}
	\textstyle\frac{\mathrm{d}E_\mathrm{b}}{\mathrm{d}s}(s) = \begin{cases}
		F_\mathrm{m}(s)/\eta^+ + F_\mathrm{b,0},& \text{if }F_\mathrm{m}(s)\geq 0,\\
		F_\mathrm{m}(s)\cdot\eta^- + F_\mathrm{b,0},& \text{if }F_\mathrm{m}(s)< 0,
	\end{cases}
	\end{equation}
	
\end{small}
\endgroup
\noindent where $F_\mathrm{d}$ is the total drag force acting on the vehicle, $\kappa$ is the curvature of the spatial trajectory along the track (known~a~priori), $\eta^+$ and $\eta^-$ are the propulsive and regenerative powertrain efficiencies respectively, and $F_\mathrm{b,0}$ represents auxiliary, stationary energy demands. $F_\mathrm{m}$ is the control input corresponding to motor power demand, $F_\mathrm{brk}$ is the braking force, and $E_\mathrm{b}$ is the state tracking the cumulative (battery) energy used throughout the lap. Our use case is an electric racing car, but the methodology also applies to fuel-limited combustion engine cars, by considering $E_\mathrm{b}$ as the cumulative fuel energy used.\\
We enforce powertrain limits on the battery as
\par\nobreak\vspace{-5pt}
\begingroup
\allowdisplaybreaks
\begin{small}
	\begin{equation}
			\label{eqn4}
	\begin{aligned}
	F_\mathrm{m}(s) &\leq F^+_\mathrm{PT}(s,E_\mathrm{kin}),\\
	F_\mathrm{m}(s) &\geq F^-_\mathrm{PT}(s,E_\mathrm{kin}),\\
\end{aligned}
	\end{equation}
	
\end{small}
\endgroup
\noindent where $F^+_\mathrm{PT}$ is the maximum power output, corresponding to full-throttle operation, and $F^-_\mathrm{PT}$ is the limit for regenerative braking force, set to zero in the case of combustion engine cars.\\
In contrast to closely-related  works~\cite{SalazarBalernaEtAl2017,SalazarBalernaEtAl2018,SalazarElbertEtAl2017,EbbesenSalazarEtAl2018,NeumannFieniEtAl2023,EshofKampenEtAl2025}, we constrain grip limits by bounding the longitudinal force rather than enforcing a maximum kinetic energy bound to represent tire grip limits,
\par\nobreak\vspace{-5pt}
\begingroup
\allowdisplaybreaks
\begin{small}
		\begin{equation}
				\label{eqn5}
		\begin{aligned}
F_\mathrm{m}(s) - F_\mathrm{brk}(s) &\leq F^+_\mathrm{grip}(s,\kappa,E_\mathrm{kin}),\\
F_\mathrm{m}(s) - F_\mathrm{brk}(s) &\geq F^-_\mathrm{grip}(s,\kappa,E_\mathrm{kin}),\\
		\end{aligned}
	\end{equation}

\end{small}
\endgroup
\noindent where $F^\pm_\mathrm{grip}$ are the upper and lower grip limits, respectively. This approach ensures we exclusively have instantaneous jumps in our kinetic costate dynamics as control-independent state limits are only reached at corner apexes. We discuss this more thoroughly and explain the benefits in Section~\ref{jumpconditionssection}. Moreover, this force-limit formulation allows for a straightforward derivation from potentially data-driven \textit{g-g-v} envelope limits~\cite{DuhrSandeepEtAl2022,EshofKampenEtAl2025b,Piccinini2024,Milliken1995}, and avoids the pre-solve required to find the maximum kinetic energy bound. Our choices for defining these limits and the drag forces can be found in the \hyperref[aap]{Appendix}.\\
The full optimal control problem with states $\textbf{x}(s) \coloneq \{E_\mathrm{kin},E_\mathrm{b}\}$ and controls $\textbf{u}(s) \coloneq \{F_\mathrm{m},F_\mathrm{brk}\}$ reads as follows:
\label{ModelSection}
\begin{prob}[Energy-limited Lap-Time Optimization]\label{Problem1}
	The minimum-lap-time control strategy is the solution of
	\newlength{\mylength}
	\settowidth{\mylength}{\quad \text{Battery/Fuel Dynamics:}} 
	
	\par\nobreak\vspace{-5pt}
	\begingroup
	\allowdisplaybreaks
	\begin{small}
		
		\begin{equation*}
			\begin{aligned}
				\min_{\textbf{u}(s),\textbf{x}(s)} \quad & \int_{s_\mathrm{0}}^{s_\mathrm{f}}\frac{1}{\sqrt{2\,E_\mathrm{kin}(s)/m}}\, \mathrm{d}s\\
				\makebox[\mylength][l]{\text{s.t.}} & \\
				\makebox[\mylength][l]{\quad \text{Vehicle Dynamics:}} & (\ref{eqn2})\\
				\makebox[\mylength][l]{\quad \text{Battery/Fuel Dynamics:}} & (\ref{eqn3})\\		
				\makebox[\mylength][l]{\quad \text{Powertrain Limits:}} & (\ref{eqn4})\\	
				\makebox[\mylength][l]{\quad \text{Grip Limits:}} & (\ref{eqn5})\\
				\makebox[\mylength][l]{\quad \text{Braking Non-negativity:}} & F_\mathrm{brk}(s) \geq 0\\
				\makebox[\mylength][l]{\quad \text{Energy Limits:}} & E_\mathrm{b}(s_\mathrm{f}) \leq \Delta E_\mathrm{b,max} +  E_\mathrm{b}(s_\mathrm{0}), \\
			\end{aligned}
		\end{equation*}
		
	\end{small}
	\endgroup
\end{prob}
\noindent where $\Delta E_\mathrm{b,max}$ is the energy budget for the lap.\\
The Hamiltonian of the optimization problem is written as
\par\nobreak\vspace{-5pt}
\begingroup
\allowdisplaybreaks
\begin{small}

	\begin{multline}
			\label{eqn8}
		H(x,u,\lambda) = \frac{1}{\sqrt{2\,E_\mathrm{kin}(s)/m}} + \lambda_\mathrm{kin}\cdot\textstyle\frac{\mathrm{d}E_\mathrm{kin}}{\mathrm{d}s}	+ \lambda_\mathrm{b}\cdot\textstyle\frac{\mathrm{d}E_\mathrm{b}}{\mathrm{d}s}\\
		+ \mu^+_\mathrm{grip}(F_\mathrm{m} - F_\mathrm{brk} -  F^+_\mathrm{grip}) - \mu^-_\mathrm{grip}(F_\mathrm{m} - F_\mathrm{brk} - F^-_\mathrm{grip})\\
		+ \mu^+_\mathrm{PT}(F_\mathrm{m}-F^+_\mathrm{PT}) - \mu^-_\mathrm{PT}(F_\mathrm{m}-F^+_\mathrm{grip}) - \mu_\mathrm{brk}F_\mathrm{brk},
	\end{multline}
	
\end{small}
\endgroup
\noindent where $\lambda_\mathrm{\{kin,b\}}$ are the costates and $\mu^\pm_\mathrm{\{grip,PT,brk\}}$ are multipliers for the respective constraints. Note that for brevity's sake, in this equation and the remainder of this paper, we do not write out the full function notations (dependencies on $s,\,E_\mathrm{kin}$). Next, in accordance with PMP, we derive our costate dynamics as the negative partial derivatives of the Hamiltonian with respect to the states:
\par\nobreak\vspace{-5pt}
\begingroup
\allowdisplaybreaks
\begin{small}
	
	\begin{equation}
			\label{eqn9}
		\begin{aligned}
			\textstyle\frac{\mathrm{d} \lambda_\mathrm{kin}}{\mathrm{d}s} = -\textstyle\frac{\partial H}{\partial E_\mathrm{kin}} &= 
			(2\,E_\mathrm{kin}/m)^{-\frac{3}{2}}/m + \lambda_\mathrm{kin}\textstyle\frac{\partial F_\mathrm{d}}{\partial E_\mathrm{kin}}
			\\
			&-\mu_\mathrm{grip}^+\textstyle\frac{\partial F_\mathrm{grip}^+}{\partial E_\mathrm{kin}} -\mu_\mathrm{grip}^-\textstyle\frac{\partial F_\mathrm{grip}^-}{\partial E_\mathrm{kin}}\\
			&-\mu_\mathrm{PT}^+\textstyle\frac{\partial F_\mathrm{PT}^+}{\partial E_\mathrm{kin}} -\mu_\mathrm{PT}^-\textstyle\frac{\partial F_\mathrm{PT}^-}{\partial E_\mathrm{kin}},\\
			\textstyle\frac{\mathrm{d} \lambda_\mathrm{b}}{\mathrm{d}s} = -\textstyle\frac{\partial H}{\partial E_\mathrm{b}} &= 0,\hfill
		\end{aligned}
	\end{equation}
	
\end{small}
\endgroup
\noindent Using the end conditions, the fact that the battery costate is constant and dual feasiblity of the constraint multipliers ($\mu\geq0$), we can derive the battery costate sign:
\par\nobreak\vspace{-5pt}
\begingroup
\allowdisplaybreaks
\begin{small}
	
	\begin{equation}
					\label{eqn10}
		\lambda_\mathrm{b}(s_\mathrm{f}) = \mu_\mathrm{b,max} \implies \lambda_\mathrm{b}(s) =  \mu_\mathrm{b,max} \geq 0 \quad \forall s.
	\end{equation}
	
\end{small}
\endgroup
\noindent Next, we define our stationarity conditions as follows:
\par\nobreak\vspace{-5pt}
\begingroup
\allowdisplaybreaks
\begin{small}
	
	\begin{equation}
					\label{eqn11}
		\begin{aligned}
			0 \in \textstyle\frac{\partial H}{\partial F_\mathrm{m}} &= \begin{cases}
				\lambda_\mathrm{kin}+\lambda_\mathrm{b}\cdot \left(\eta^+ + F_\mathrm{m}\frac{\partial\eta^+}{\partial F_\mathrm{m}}\right)/\eta^{+2}\\
				+ \mu_\mathrm{grip}^+ - \mu_\mathrm{grip}^- + \mu_\mathrm{PT}^+ - \mu_\mathrm{PT}^-,
				& \text{if } F_\mathrm{m}\geq 0,\\
				\lambda_\mathrm{kin}+\lambda_\mathrm{b}\cdot\left(\eta^- + F_\mathrm{m}\frac{\partial\eta^-}{\partial F_\mathrm{m}}\right)\\
				+ \mu_\mathrm{grip}^+ - \mu_\mathrm{grip}^- + \mu_\mathrm{PT}^+ - \mu_\mathrm{PT}^-, & \text{if } F_\mathrm{m}< 0,\\
			\end{cases}\\
			0 \in \textstyle\frac{\partial H}{\partial F_\mathrm{brk}} &= -\lambda_\mathrm{kin} - \mu_\mathrm{grip}^+ + \mu_\mathrm{grip}^- -\mu_\mathrm{brk}.
		\end{aligned}
	\end{equation}
	
\end{small}
\endgroup
\noindent Complementary slackness dictates that:
\par\nobreak\vspace{-5pt}
\begingroup
\allowdisplaybreaks
\begin{small}
	
	\begin{equation}
					\label{eqn12}
		\begin{aligned}
			\mu^+_\mathrm{grip}\left(F_\mathrm{m}-F_\mathrm{brk}-F^+_\mathrm{grip}\right) &= 0,\\
						\mu^-_\mathrm{grip}\left(-F_\mathrm{m}+F_\mathrm{brk}+F^+_\mathrm{grip}\right) &= 0,\\
									\mu^+_\mathrm{PT}\left(F_\mathrm{m}-F^+_\mathrm{PT}\right) &= 0,\\
						\mu^-_\mathrm{PT}\left(-F_\mathrm{m}+F^+_\mathrm{PT}\right) &= 0,\\								
						\mu_\mathrm{brk}\left(-F_\mathrm{brk}\right) &= 0.\\							
		\end{aligned}	
	\end{equation}
	
\end{small}
\endgroup

\noindent From these two conditions follows a policy for control inputs as a function of costate variables. In motorsport applications, the vehicle operates most of the lap at the fixed operating points of the powertrain limits, corresponding to full-throttle operation or maximum regenerative braking. Therefore, we can assume $\textstyle\frac{\partial \eta^\pm}{\partial F_\mathrm{m}}=0$, i.e., that the powertrain efficiency is independent of the force demands, by linearizing the loss model around these operating points. Through this assumption, we obtain the following optimal control policy as a function of the costate variables:
\par\nobreak\vspace{-5pt}
\begingroup
\allowdisplaybreaks
\begin{small}
	
	\begin{equation}
					\label{eqn13}
	\begin{aligned}
		F_\mathrm{m} &\begin{cases}
			= \min(F^+_\mathrm{PT},F^+_\mathrm{grip}), & \text{if } \lambda_\mathrm{kin}/\lambda_\mathrm{b}	< -1/\eta^+,\\
			\in [0,\min(F^+_\mathrm{PT},F^+_\mathrm{grip})], & \text{if } \lambda_\mathrm{kin}/\lambda_\mathrm{b} = -1/\eta^+,\\
			= 0, & \text{if } -1/\eta^+ < \lambda_\mathrm{kin}/\lambda_\mathrm{b}	< -\eta^-,\\
						\in [\max(F^-_\mathrm{PT},F^-_\mathrm{grip}),0], & \text{if } \lambda_\mathrm{kin}/\lambda_\mathrm{b} = -\eta^-,\\
						= \max(F^-_\mathrm{PT},F^-_\mathrm{grip}), & \text{if } -\eta^-<\lambda_\mathrm{kin}/\lambda_\mathrm{b},\\
		\end{cases}\\
		F_\mathrm{brk} &\begin{cases}
	 	= 0, & \text{if } \lambda_\mathrm{kin}	< 0,\\ 		
	 	= \max(F^-_\mathrm{PT},F^-_\mathrm{grip})-F^-_\mathrm{grip}, & \text{if } \lambda_\mathrm{kin}	> 0,\\ 		
	 	\in [0,\max(F^-_\mathrm{PT},F^-_\mathrm{grip})-F^-_\mathrm{grip}], & \text{if } \lambda_\mathrm{kin} = 0.\\ 		
		\end{cases}
	\end{aligned}
	\end{equation}
	
\end{small}
\endgroup

\begin{figure}[t]
	\centering
	\includegraphics[width=\linewidth]{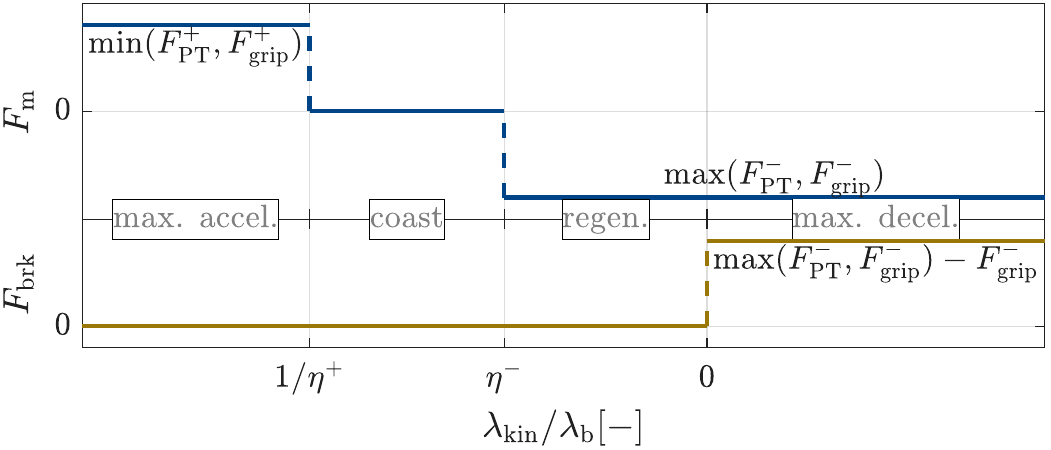}
	\caption{Optimal bang-bang control policy as a function of the costate ratio.}
	\label{fig:policy}
\end{figure}
\noindent This policy is visualized in Fig.~\ref{fig:policy} and shows that by selecting linear loss models for the powertrain, the optimal control policy follows a bang-bang strategy with distinct phases: When $\lambda_\mathrm{kin}/\lambda_\mathrm{b} < -1/\eta^+$, the driver is responsible for maximizing acceleration, following the grip limit by modulating the throttle or activating the powertrain limits by fully depressing the throttle pedal. When $-1/\eta^+<\lambda_\mathrm{kin}/\lambda_\mathrm{b}<-\eta^-$, the driver should coast, i.e. lift off the throttle without any braking action. When $-\eta^-<\lambda_\mathrm{kin}/\lambda_\mathrm{b}$, the driver should maximize regenerative braking, and avoid using the mechanical friction brakes until the latest possible braking point, i.e., where the grip limit activates again, to decelerate and complete the corner. These distinct phases allow for safe usage alongside human drivers, as discrete signals allow the driver to maintain control without requiring an active overwrite of throttle pedal or brake inputs.\\
To finalize our formulations, we derive the equations for our constraint multipliers using the control policy (6) and stationarity conditions (4) as
\par\nobreak\vspace{-5pt}
\begingroup
\allowdisplaybreaks
\begin{small}
	
	\begin{equation}
					\label{eqn14}
		\begin{aligned}
			\mu_\mathrm{grip}^+ &= \begin{cases}
				-(\lambda_\mathrm{kin}+\lambda_\mathrm{b}/\eta^+ ), & \text{if } \lambda_\mathrm{kin}/\lambda_\mathrm{b}	< -1/\eta^+\\
				& \text{and } F_\mathrm{grip}^+ < F_\mathrm{PT}^+\\
				0, & \text{otherwise},
			\end{cases}\\
			\mu_\mathrm{grip}^- &= \begin{cases}
	\lambda_\mathrm{kin}+\lambda_\mathrm{b}\cdot\eta^-, & \text{if } -\eta^-<\lambda_\mathrm{kin}/\lambda_\mathrm{b} \\
	& \text{and } F_\mathrm{grip}^- > F_\mathrm{PT}^-,\\
		\lambda_\mathrm{kin}, & \text{if } 0<\lambda_\mathrm{kin}/\lambda_\mathrm{b} \\
	& \text{and } F_\mathrm{grip}^- < F_\mathrm{PT}^-,\\
0, & \text{otherwise},
\end{cases}\\
			\mu_\mathrm{PT}^+ &= \begin{cases}
	-(\lambda_\mathrm{kin}+\lambda_\mathrm{b}/\eta^+ ), & \text{if } \lambda_\mathrm{kin}/\lambda_\mathrm{b}	< -1/\eta^+\\
	& \text{and } F_\mathrm{PT}^+ < F_\mathrm{grip}^+,\\
0, & \text{otherwise},
\end{cases}\\
\mu_\mathrm{PT}^- &= \begin{cases}
-(\lambda_\mathrm{kin}+\lambda_\mathrm{b}\cdot\eta^- ), & \text{if } -\eta^-<\lambda_\mathrm{kin}/\lambda_\mathrm{b}<0 \\
& \text{and } F_\mathrm{PT}^- > F_\mathrm{grip}^-,\\
\lambda_\mathrm{b}\cdot \eta^-,& \text{if } 0<\lambda_\mathrm{kin}/\lambda_\mathrm{b}\\
& \text{and } F_\mathrm{PT}^- > F_\mathrm{grip}^-,\\
0, & \text{otherwise},
\end{cases}\\
\mu_\mathrm{brk} &= \begin{cases}
	0, & \text{if } 0<\lambda_\mathrm{kin}/\lambda_\mathrm{b} \\
	& \text{and } F_\mathrm{PT}^- > F_\mathrm{grip}^-,\\
	\lambda_\mathrm{b}/\eta^+, &  \text{if } \lambda_\mathrm{kin}/\lambda_\mathrm{b}	< -1/\eta^+\\
	& \text{and } F_\mathrm{grip}^+ < F_\mathrm{PT}^+,\\
	\lambda_\mathrm{b}\cdot\eta^-, &  \text{if } \eta^-<\lambda_\mathrm{kin}/\lambda_\mathrm{b}	< 0\\
	&  \text{or } (0<\lambda_\mathrm{kin}/\lambda_\mathrm{b} \text{ and } F_\mathrm{grip}^- > F_\mathrm{PT}^-),\\
	-\lambda_\mathrm{kin}, & \text{otherwise}.
\end{cases}\\		\end{aligned}
	\end{equation}
	
\end{small}
\endgroup
\subsection{Singularities}
We note that our control policy in (\ref{eqn13}) is not explicitly defined at the distinct switching points, $\lambda_\mathrm{kin}=\{-\lambda_\mathrm{b}/\eta^+,-\lambda_\mathrm{b}\cdot\eta^-,0\}$. In practice, this indeterminacy occurs only if there are costate equilibria at these points. To form an equilibrium, $\textstyle\frac{\mathrm{d} \lambda_\mathrm{kin}}{\mathrm{d}s}=0$ must hold, and solving this for $E_\mathrm{kin}$ gives the singular speeds where our policy is not yet explicitly defined:
\par\nobreak\vspace{-5pt}
\begingroup
\allowdisplaybreaks
\begin{small}
	
	\begin{equation}
					\label{eqn15}
		E_\mathrm{kin,sing} = \frac{m}{2}\left(-\frac{\partial F_\mathrm{d}}{\partial E_\mathrm{kin}} m \{-\lambda_\mathrm{b}/\eta^+,-\lambda_\mathrm{b}\cdot\eta^-,0\}\right)^{-\frac{2}{3}}.
	\end{equation}
	
\end{small}
\endgroup
\noindent In these costate equilibria, the optimal solution attempts to track this critical kinetic energy. This corresponds to the undesirable situation where, due to extreme energy limitations, a pure coast-regen strategy is not feasible and the driver is required to modulate throttle application to optimally minimize lap-time. In practice, we avoid this situation either by driving fewer laps and allowing more energy per lap, or by tightening the upper powertrain bound $F_\mathrm{PT}^+$, i.e. reducing full-throttle power output, such that the critical velocity is not met. In our implementation we therefore use definition (\ref{eqn15}) to detect infeasible energy limitations, and return a warning to either increase energy budget or reduce full-throttle power output. Alternatively, the car could switch to a speed control mode and override driver power similar to how a road-car cruise control system would, which could still be a safe option if changes to $E_\mathrm{kin,sing}$ are minor.
In any case, this consideration completes our definition of the continuous optimal solution dynamics.
\subsection{Jump Conditions}
\label{jumpconditionssection}
Equations (\ref{eqn2}), (\ref{eqn3}), (\ref{eqn9}), (\ref{eqn13}), and (\ref{eqn14}) fully define the continuous optimal dynamics. For tire grip there are, however, combinations of kinetic energy and curvature where there exist no feasible values for the longitudinal tire forces $F_\mathrm{m}$ and $F_\mathrm{brk}$. This is evident from our definitions found in the \hyperref[aap]{Appendix}, where we also derive a formulation for the kinetic energy where these limits are reached. These limits activate at speeds where all available tire grip is devoted to lateral forces, corresponding to $F_\mathrm{grip}^+=F_\mathrm{grip}^-=0$, which limits the car's maximum speed through corners. As these are inherent upper limitations to the kinetic energy state, activating this constraint results in a discontinuous jump in the kinetic costate value~\cite{KeulenGillotEtAl2014}. Accordingly, to ensure the longitudinal tire forces are at zero at the jump point, the jump conditions are as follows:
\par
\begin{figure}[t]
	\centering
	\includegraphics[width=\linewidth]{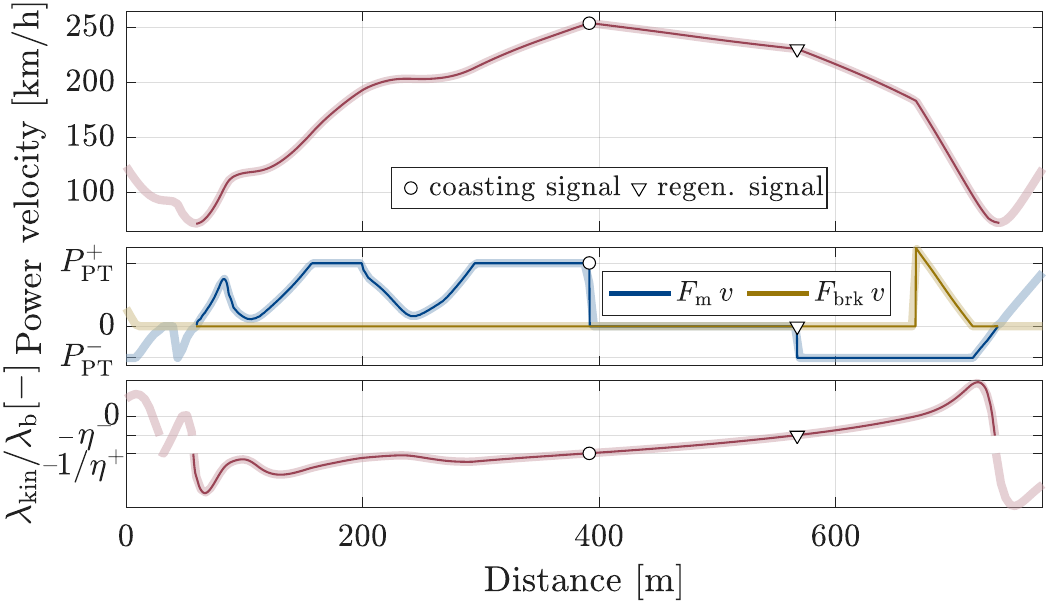}
	\caption{Optimal trajectories for a single straight, between two costate jumping points at the apexes of subsequent corners, found by integrating the optimal solution dynamics. Two discrete signals suffice to instruct the driver to manage energy, minimizing cognitive load and distraction. The result of a direct solving method (NLP) on the problem is shown in a lighter color, matching and verifying the derivation of the costate dynamics.}
	\label{fig:singlestraight}
\end{figure}

\nobreak\vspace{-5pt}
\begingroup
\allowdisplaybreaks
\begin{small}
	
	\begin{equation}
		\begin{aligned}
			\lambda_\mathrm{kin}(s_i^-) &\begin{cases}
				\in [\lambda_\mathrm{kin}(s_i^+),-\lambda_\mathrm{b}\cdot \eta^-], & \text{if } F_\mathrm{m}(s^-)-F_\mathrm{brk}(s^-) = 0,\\
				= -\lambda_\mathrm{b}\cdot \eta^-, & \text{if } F_\mathrm{m}(s^-)-F_\mathrm{brk}(s^-) < 0,\\
			\end{cases}\\
			\lambda_\mathrm{kin}(s_i^+) &\begin{cases}
				\in [-\lambda_\mathrm{b}/ \eta^+,\lambda_\mathrm{kin}(s_i^-)], & \text{if } F_\mathrm{m}(s^+)-F_\mathrm{brk}(s^+) = 0,\\
				= -\lambda_\mathrm{b}/ \eta^+, & \text{if } F_\mathrm{m}(s^+)-F_\mathrm{brk}(s^+) > 0,\\
			\end{cases}
		\end{aligned}
	\end{equation}
	
\end{small}
\endgroup
\noindent where $s_i^\pm$ are the distances infinitely close before and after the jumping point. Using our continuous dynamics, we can simulate and predict the optimal solution between two subsequent jumping points, as visualized in Fig.~\ref{fig:singlestraight}.
\renewcommand{\algorithmicendwhile}{\textbf{end}}
\begin{algorithm}[t]
	\caption{Solving procedure of Problem~\ref{Problem1} given $\lambda_\mathrm{b}$}\label{alg1}
	\begin{algorithmic}
	\STATE \text{ignore energy limits}\\
	\STATE \text{set $s_0$ to highest curvature corner apex}\\
	\STATE $E_\mathrm{kin}(s_0) = \text{arg min}_{E_\mathrm{kin}}|F_\mathrm{grip}^+(s_0)-F_\mathrm{grip}^-(s_0)|$
	\STATE $i\gets 0$
	\WHILE{$s_i \leq s_\mathrm{f}$} 
	\STATE \text{minimize $\lambda_\mathrm{kin}(s_i^+)$ by repeated integration of $s \in [s_i,s_\mathrm{f}]$}\\
	\STATE \text{until apex grip limits are violated} \COMMENT{bisection method}\\
	\STATE $s_{i+1} = \text{arg min}_{s \in [s_i s_\mathrm{f}]}|F_\mathrm{grip}^+-F_\mathrm{grip}^-|$\COMMENT{next corner apex}\\
	\text{append $\{x^\star,u^\star,\lambda_\mathrm{kin}^\star | s\in[s_i,s_{i+1}]\}$ to optimal solution}\\
	\STATE $i\gets i+1$
	\ENDWHILE
	\end{algorithmic}
\end{algorithm}%

\subsection{Solution Algorithm}
By deriving equations for our continuous optimal solution dynamics, we have reduced the optimal control problem to finding the value of our energy limitation, i.e., the energy limit constraint multiplier $\lambda_\mathrm{b}=\mu_\mathrm{b,max}$, and the magnitude of the instantaneous jumps in kinetic costate at the small number of corner apexes, significantly reducing computational complexity. This reduces the problem to a set of sequential boundary value problems, which we can solve with single shooting methods~\cite{KeulenGillotEtAl2014}. Moreover, we can exploit monotonicity and use efficient root-finding methods to quickly derive the optimal solution. Increasing the magnitude of the jump leads to later coasting and later braking, and vice versa. Intuitively, this is analogous to the driver braking too early or too late for the corner, meaning that by checking whether the car violated maximum apex speed, we know whether to increase of reduce the magnitude of the jump. Our proposed solving method is outlined in Algorithm~\ref{alg1}, using a bisection method for minimization of the kinetic energy costate value after the jump. Algorithm~\ref{alg1} can be nested within another bisection algorithm that maximizes $\lambda_\mathrm{b}$ until the $E_\mathrm{b}(s_\mathrm{f}) \leq \Delta E_\mathrm{b,max}+E_\mathrm{b}(s_\mathrm{0})$ constraint is violated~\cite{KeulenGillotEtAl2014}, thereby identifying $\lambda_\mathrm{b}$ and fully solving the problem.
\section{Results}
\label{ResultsSection}
\begin{figure}[t]
	\centering
	\includegraphics[width=\linewidth]{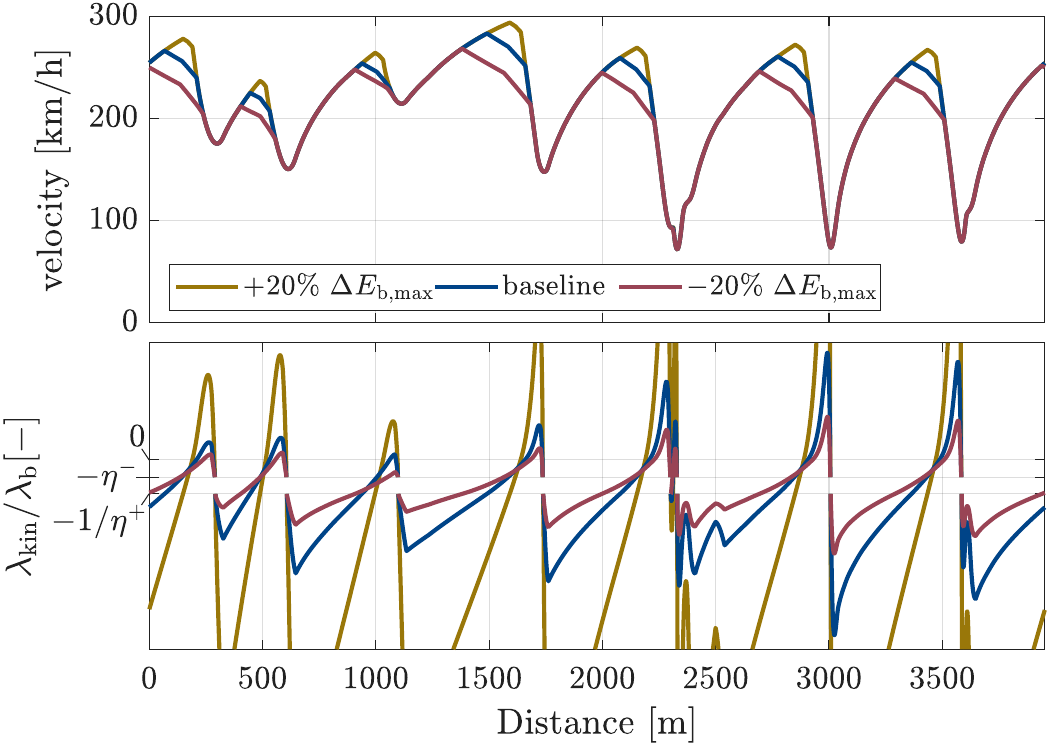}
	\caption{Velocity and costate ratio trajectories over a full lap of the Zolder circuit for three different energy budgets. Lap-time-optimal coasting and regenerative braking phases are triggered as the costate ratio passes the thresholds.}
	\label{fig:fulllap}
\end{figure}

In this section we evaluate the performance of our proposed method. We compare our method to a numerical solution obtained through a nonlinear programming solver, in terms of solving times and accuracy. Next, we explore the potential advantages of an adaptive kinetic costate trajectory for large differences in energy limitations. First, we look at the results for a full lap of the Zolder circuit for three different energy budgets, as visualized in Fig.~\ref{fig:fulllap}. Dictated by the costate ratio trajectory, optimal coasting and regenerative braking phases are identified and, for the full lap, only seven instantaneous coasting cues, followed by regeneration cues, are required to be passed to the driver per lap for lap-time-optimal energy management. Another visualization of the optimal coast and regeneration points is provided in Fig.~\ref{fig:track}, showing that the optimal coasting and regeneration points are always near the ends of straight sections, where lap-time sensitivity to varying energy deployment is small. This means that the signals are passed to the driver when the car is vehicle-dynamically stable and the driver is not preoccupied with executing a corner at the limits of the car, where such a signal may be an undesirable distraction.
\subsection{Comparison to Direct Methods}
Previous works solved the energy management optimal control problem using direct methods like nonlinear programming. As shown in Fig.~\ref{fig:singlestraight}, the same solution is obtained through our indirect method, verifying our derivations of the costate dynamics. The advantage of our indirect solution methodology is twofold: First, it avoids reliance on third-party solvers. Second, the simplicity of the solving algorithm makes it lightweight and very efficient, bringing the solving time down from the order of seconds~\cite{KampenHerrmannEtAl2023,EshofKampenEtAl2025} to the order of milliseconds, paving the way for real-time adaptations to changing conditions during a race.
\subsection{Advantages of an Optimal Costate Trajectory}
\begin{figure}[t]
	\centering
	\includegraphics[width=\linewidth]{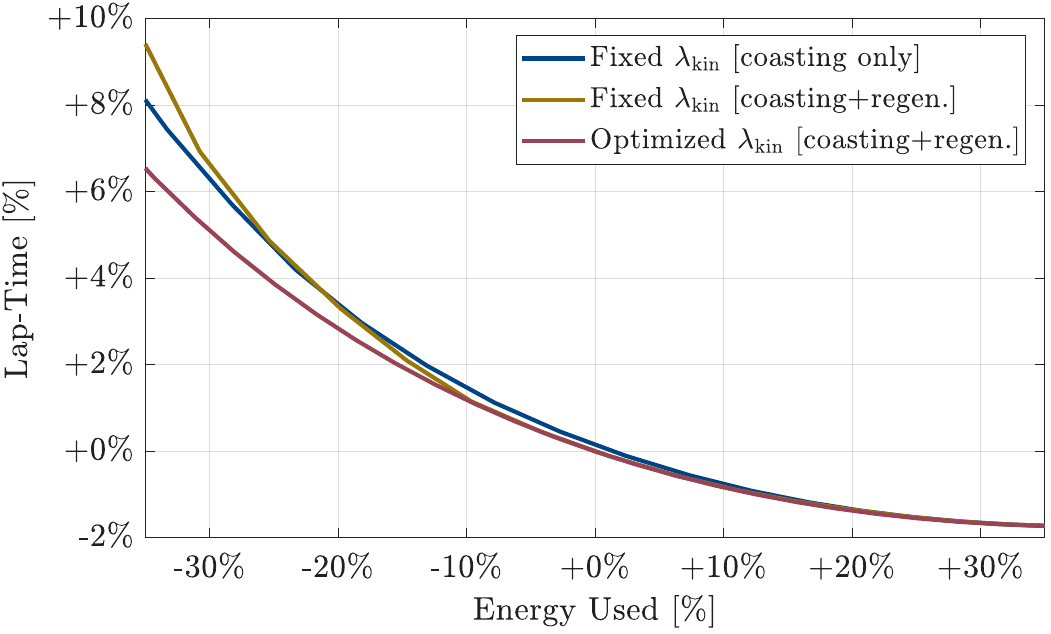}
	\caption{Lap-time gains over energy budget for the different approaches. An adaptive costate trajectory is very advantageous, and its benefits are bigger when energy limits are tighter. For fixed costate trajectories and severe energy limitations, a regenerative braking phase may be disadvantageous due to the car slowing down too quickly and early in some corners.}
	\label{fig:comparison}
\end{figure}

In previous implementations on cars~\cite{SalazarBalernaEtAl2018,EshofKampenEtAl2025}, the costate trajectory was assumed fixed for a simple and lightweight integration, and to avoid reliance on third party solvers. In these implementations, model errors and other disturbances are dealt with by adjusting the costate thresholds using a feedback loop. As we demonstrate in this paper, the costate dynamics are influenced by when and where constraints are activated and depend on $\lambda_\mathrm{b}$, i.e., the severity of the energy limitation. Accordingly, assuming an unchanging $\lambda_\mathrm{kin}$ trajectory for a changing $\lambda_\mathrm{b}$ leads to sub-optimal results. This is visualized in Fig.~\ref{fig:comparison}, where we compare our optimized $\lambda_\mathrm{kin}$ method to one using a fixed $\lambda_\mathrm{kin}$ trajectory based on the nominal case, but with adapted thresholds, and the method from~\cite{EshofKampenEtAl2025}, which only considers a coasting phase. The optimized $\lambda_\mathrm{kin}$ method is globally optimal and serves as a benchmark of the other two methods. The method that exclusively considers a coasting phase is outperformed across the range of energy budgets, but does outperform the method that also includes a regeneration phase for stronger energy limitations. This is because the fixed costate trajectory may lead to premature regeneration signals, causing the driver to slow down too quickly and early for a corner, resulting in unnecessary extra lap-time.
\section{Conclusion} \label{Conclusion}
This paper presented a provably optimal, real-time capable method for energy management of racing cars. By deriving the optimal control policy and costate dynamics, we turned the energy management problem into a set of boundary value problems that we could efficiently solve using simple algorithms. Moreover, we explicitly determined the necessary conditions for the optimal control solution to follow a human-executable bang-bang trajectory, with control cues activating as the costate trajectory passes powertrain dependent thresholds. We compared our solution method to a nonlinear programming solution, showing that, whilst the solutions are identical, our proposed method solves orders of magnitude faster, paving the way for adaptive on-board implementations. Moreover, we showed that an optimal kinetic costate trajectory can significantly improve lap-time compared to previously proposed fixed-costate approaches.\\
In future work, we plan to validate human-in-the-loop implementation of our method through racetrack testing, whereby adaptive control schemes could be investigated to account for disturbances and changing conditions.
Moreover, the approach is not limited to our fully-electric case-study, but can be readily extended and applied to hybrid power units.


\appendix
\label{aap}
\section*{Vehicle Dynamic Model and Constraints}
This section includes the vehicle dynamic model parametrizations and formulations we used throughout this paper. More generally, the proposed methodology could also be applied to other formulations.\\
First, we define the total drag force acting on the vehicle as
\par\nobreak\vspace{-5pt}
\begingroup
\allowdisplaybreaks
\begin{small}
	
	\begin{multline}
		 F_\mathrm{d}(s,\kappa,E_\mathrm{kin}) = c_\mathrm{d}\,A\,\rho\,E_\mathrm{kin}(s)/m + c_\mathrm{\kappa}|\kappa(s)|E_\mathrm{kin}(s)\\ + c_\mathrm{r}\,F_\mathrm{z}(s,E_\mathrm{kin}),
	\end{multline}
	
\end{small}
\endgroup
\noindent where $c_\mathrm{d}$ describes air drag coefficient, $A$ is aerodynamic reference area, $\rho$ is the air density, $c_\mathrm{\kappa}$ is a coefficient describing cornering resistances, $c_\mathrm{r}$ describes rolling resistances, and $F_\mathrm{z}$ is the vertical force\cite{EbbesenSalazarEtAl2018,GuzzellaSciarretta2007}. The vertical force is defined as
\par\nobreak\vspace{-5pt}
\begingroup
\allowdisplaybreaks
\begin{small}
	\begin{equation}
		F_\mathrm{z}(s,E_\mathrm{kin}) = c_\mathrm{l}\,A\,\rho\,E_\mathrm{kin}(s)/m + m\,g,
	\end{equation}
\end{small}%
\endgroup
\noindent where $c_\mathrm{l}$ is the downforce coefficient and $g$ is the gravitational acceleration.\\
The powertrain limits are defined as power limits:
\par\nobreak\vspace{-5pt}
\begingroup
\allowdisplaybreaks
\begin{small}
	\begin{equation}
		F_\mathrm{PT}^{\pm}(s,E_\mathrm{kin}) = P^{\pm}_\mathrm{PT} / \sqrt{2\,E_\mathrm{kin}(s)/m},
	\end{equation}
\end{small}%
\endgroup
\noindent where $ P^{\pm}_\mathrm{PT}$ are upper and lower motor power limits, respectively.\\
Finally, the grip limits are defined as friction ellipses, accounting for downforce:
\par\nobreak\vspace{-5pt}
\begingroup
\allowdisplaybreaks
\begin{small}
	\begin{equation}
		F_\mathrm{grip}^\pm(s,\kappa,E_\mathrm{kin}) = \pm \mu_\mathrm{x} \sqrt{F_\mathrm{z}(s,E_\mathrm{kin})^2 - (F_\mathrm{y}(s,\kappa,E_\mathrm{kin})/\mu_\mathrm{y})^2},
	\end{equation}
\end{small}%
\endgroup
\noindent where $\mu_\mathrm{\{x,y\}}$ are longitudinal and lateral grip coefficient, respectively, and $F_\mathrm{y}$ is the lateral (centrifugal) force acting on the vehicle, which is found through the curvature and kinetic energy as
\par\nobreak\vspace{-5pt}
\begingroup
\allowdisplaybreaks
\begin{small}
	
	\begin{equation}
		F_\mathrm{y}(s,\kappa,E_\mathrm{kin}) = 2\,\kappa(s)\,E_\mathrm{kin}(s).
	\end{equation}
	
\end{small}%
\endgroup
\noindent Through these formulations, we can also derive the absolute highest kinetic energy possible without violating grip limits, i.e., the kinetic energy limits where costate jumps occur:
\par\nobreak\vspace{-5pt}
\begingroup
\allowdisplaybreaks
\begin{small}
	
	\begin{equation}
E_\mathrm{kin,jump}(s,\kappa) =
\begin{cases}
	\dfrac{mg}{\dfrac{2|\kappa(s)|}{\mu_y} - \dfrac{c_l A \rho}{m}}, 
	& \text{if } \dfrac{2|\kappa(s)|}{\mu_y} - \dfrac{c_l A \rho}{m} > 0, \\[1.2em]
	+\infty, 
	& \text{otherwise}.
\end{cases}
	\end{equation}
	
\end{small}%
\endgroup
\newpage
\bibliographystyle{IEEEtran}
\renewcommand{\baselinestretch}{0.96}
\bibliography{../../../bibliography/powertrains,../../../bibliography/main,../../../bibliography/SML_papers}

\end{document}